\documentclass[MSNbibl,citesort,nameyear,dvips]{arxbj}
\usepackage{graphicx}
\aid{0}
\volume{19}
\issue{4}
\pubyear{2013}
\firstpage{1484}
\lastpage{1500}
\doi{10.3150/13-BEJSP14} 

\makeatletter
\newproclaim{result}{Result}
\makeatother

\begin{document}
\begin{frontmatter}

\title{Stability}
\runtitle{Stability}

\begin{aug}
\author{\fnms{Bin} \snm{Yu}\corref{}\ead[label=e1]{binyu@stat.berkeley.edu}}
\runauthor{B. Yu} 
\address{Departments of Statistics and EECS, University of California at
Berkeley, Berkeley, CA 94720, USA. \printead{e1}}
\end{aug}


%
\begin{abstract}
Reproducibility is imperative for any scientific discovery. More often
than not,
modern scientific findings rely on statistical analysis of
high-dimensional data. At
a minimum, reproducibility manifests itself in stability of statistical results
relative to ``reasonable'' perturbations to data and to the model used.
Jacknife, bootstrap, and cross-validation are based on perturbations to
data, while robust statistics methods deal with perturbations to models.

In this article, a case is made for the importance of stability in
statistics. Firstly, we motivate the necessity of stability for
interpretable and reliable encoding models
from brain fMRI signals. Secondly,
we find strong evidence in the literature to demonstrate
the central role of stability in statistical inference, such as
sensitivity analysis and effect detection.
Thirdly, a smoothing parameter selector based on estimation stability (ES),
ES-CV, is proposed for Lasso, in order to
bring stability to bear on cross-validation (CV). ES-CV is then
utilized in the encoding
models to reduce the number of predictors by 60\% with almost no
loss (1.3\%) of prediction performance across over 2,000 voxels.
Last, a novel ``stability'' argument is seen to drive
new results that shed light on the intriguing interactions
between sample to sample variability and
heavier tail error distribution (e.g., double-exponential)
in high-dimensional regression models with $p$ predictors and $n$ independent
samples.
In particular, when $p/n \rightarrow\kappa\in(0.3,1)$ and
the error distribution is double-exponential, the Ordinary Least
Squares (OLS) is a better estimator than
the Least Absolute Deviation (LAD) estimator.
\end{abstract}

%

\end{frontmatter}

\section{Introduction}
In his seminal paper ``The Future of Data Analysis'' (Tukey, \citeyear{Tuk62}), John
W. Tukey writes:

\begin{quote}
``It will still be true that there will be aspects of data analysis
well called technology, but there will also be the hallmarks of
stimulating science: intellectual adventure, demanding calls upon
insight, and a need to find out `how things really are' by
investigation and the confrontation of insights with experience'' (p.~63).
\end{quote}

Fast forward to 2013 in the age of information technology,
these words of Tukey ring as true as fifty years ago, but
with a new twist: the ubiquitous and massive data today were
impossible to imagine in 1962.
From the point of view of science,
information technology and data are a blessing, and a curse.
The reasons for them to be a blessing are
many and obvious. The reasons for it to be a curse
are less obvious. One of them is well articulated recently by
two prominent biologists in an editorial Casadevall and Fang (\citeyear{CasFan11})
in \emph{Infection and Immunity} (of the American Society
for Microbiology):

\begin{quote}
``Although scientists have always comforted themselves with the thought
that science is self-correcting, the
immediacy and rapidity with which knowledge disseminates today means
that incorrect information
can have a profound impact before any corrective process can take
place'' (p. 893).
\end{quote}


\begin{quote}
``A recent study analyzed the cause of retraction for 788 retracted
papers and found that error and fraud were responsible for 545 (69\%)
and 197 (25\%) cases, respectively, while the cause was unknown in 46
(5.8\%) cases~(31)'' (p. 893).
\end{quote}

The study referred is Steen (\citeyear{Ste11a}) in the \emph{Journal
of Medical Ethics}. Of the 788 retracted papers from PubMed from 2000
to 2010,
69\% are marked as ``errors'' on the retraction records.
Statistical analyses are likely to be involved in these errors.
Casadevall and Fang go on to call for ``enhanced training in probability
and statistics,'' among other remedies including ``reembracing philosophy.''
More often than not, modern scientific findings rely on statistical
analyses of high-dimensional data, and reproducibility is imperative
for any scientific discovery.
Scientific reproducibility therefore is a responsibility
of statisticians.
At a minimum, reproducibility manifests itself in the stability of
statistical results relative to ``reasonable'' perturbations to data and
to the method or model used.

Reproducibility of scientific conclusions is closely related
to their reliability. It is receiving much well-deserved attention
lately in the
scientific community (e.g., Ioannidis, \citeyear{Ioa05}; Kraft et al., \citeyear{KraZegIoa09},
Casadevall and Fang, \citeyear{CasFan11}; Nosek et al., \citeyear{nosek12}) and in the media (e.g.,
Naik, \citeyear{Naik11}; Booth, \citeyear{Boo12}).
Drawing a scientific conclusion involves multiple steps.
First, data are collected by one laboratory or one group, ideally with
a clear
hypothesis in the mind of the experimenter or scientist. In the age
of information technology, however, more and more massive amounts of data
are collected for fishing expeditions to ``discover'' scientific facts.
These expeditions involve running computer codes on data
for data cleaning and analysis (modeling and validation).
Before these facts become ``knowledge,''
they have to be reproduced or replicated
through new sets of data by the same group or preferably by other groups.
Given a fixed set of data, Donoho et al. (\citeyear{DonMalSha09}) discuss reproducible research
in computational
hormonic analysis with implications on computer-code or computing-environment
reproducibility in computational sciences including statistics.
Fonio et al. (\citeyear{FonGalBen12}) discuss replicability between laboratories
as an important screening mechanism for discoveries.
Reproducibility could have multitudes of meaning to different people.
One articulation on the meanings
of reproducibility, replication, and repeatability can be found in Stodden
(\citeyear{Stod11}).

In this paper, we advocate for more involvement of statisticians in
science and for an enhanced emphasis on
stability within the statistical framework.
Stability has been of a great concern in statistics.
For example, in the words of Hampel et al. (\citeyear{HamRonRou86}), ``\ldots robustness
theories can be viewed as stability theories of statistical inference''
(p. 8).
Even in low-dimensional linear regression models,
collinearity is known to cause instability of OLS
or problem for individual parameter estimates so that significance
testing for these estimates becomes unreliable.
Here we demonstrate the importance of statistics for understanding our brain;
we describe our methodological work on estimation stability that helps
interpret models reliably in neuroscience; and we articulate
how our solving neuroscience problems motivates theoretical work
on stability and robust statistics in high-dimensional
regression models. In other words, we tell an interwinding story of
scientific investigation and statistical developments.

The rest of the paper is organized as follows.
In Section~\ref{Sec:fmri},
we cover our ``intellectual adventure'' into neursocience,
in collaboration with the Gallant Lab at UC Berkeley,
to understand human visual pathway via fMRI
brain signals invoked by natural stimuli (images or movies)
(cf. Kay et al., \citeyear{KayNasPre08}, Naselaris et al., \citeyear{NasPreKay09}, Kay and Gallant, \citeyear{KayGal09},
Naselaris et al., \citeyear{NasKayGal11}).
In particular, we describe
how our statistical encoding and decoding models
are the backbones of ``mind-reading computers,'' as one of
the 50 best inventions of 2011 by the Time Magazine (Nishimoto et al., \citeyear{NisVuNas11}).
In order to
find out ``how things really are,'' we argue that
reliable interpretation needs stability. We define stability
relative to a data perturbation scheme.
In Section~\ref{Sec:stability},
we briefly review the vast literature on different
data perturbation schemes such as jacknife, subsampling, and bootstrap.
(We note that data perturbation in general means not only taking subsets
of data units from a given data set, but also sampling from an
underlying distribution or replicating the experiment for a new set of data.)

In Section~\ref{Sec:ES-CV}, we review an estimation stability (ES) measure
taken from Lim and Yu (\citeyear{LimYu13})
for regression feature selection. Combining ES with CV as in Lim and Yu (\citeyear{LimYu13})
gives rise to a smoothing parameter selector ES-CV for Lasso
(or other regularization methods).
When we apply ES-CV to the movie-fMRI data, we obtain a 60\% reduction
of the model size or the number of features selected at a negligible
loss of 1.3\% in terms of prediction accuracy.
Subsequently, the ES-CV-Lasso models are both sparse and more reliable
hence better suited for interpretation due to their stability and simplicity.
The stability considerations in our neuroscience endeavors have
prompted us
to connect with the concept of stability from the robust statistics
point of view.
In El Karoui et al. (\citeyear{ElKarouietal12}), we obtain very interesting theoretical results
in high-dimensional regression models with $p$ predictors and $n$
samples, shedding light on how
sample variability in the design matrix meets heavier tail error distributions
when $p/n$ is approximately a constant in $(0,1)$ or in the random matrix
regime.
We describe these results in an important special case in Section \ref
{Sec:robust}. In particular, we see that when $p/n \rightarrow\kappa$ and
$ 1 > \kappa>0.3$ or so, the Ordinary Least Squares (OLS) estimator is better
than the Least Absolute Deviation (LAD) estimator
when the error distribution is double exponential.
We conclude in Section~\ref{Sec:conclusions}.

\section{Stable models are necessary for understanding visual pathway}
\label{Sec:fmri}

Neuroscience holds the key to understanding how our mind works.
Modern neuroscience is invigorated by massive and multi-modal
forms of data enabled by advances in technology (cf. Atkil, Martone and
Van Essen, \citeyear{AkiMarEss12}).\vadjust{\goodbreak}
Building mathematical/statistical models on this data,
computational neuroscience is at the frontier of neuroscience.
The Gallant Lab at UC Berkeley is a leading neuroscience lab
specializing in understanding the visual pathway, and is a long-term
collaborator with the author's research group. It pioneered the use of
natural stimuli
in experiments to invoke brain signals, in contrast to synthetic
signals such as white noise and moving bars or checker boards as
previously done.

Simply put, the human visual pathway works as follows.
Visual signals are recorded by retina and through the relay center LGN
they are transmitted to primary visual cortex areas V1, on to V2 and V4,
on the ``what'' pathway (in contrast to the ``where'' pathway) (cf.
Goodale and Milner, \citeyear{GooMil92}).
Computational vision neuroscience aims at modeling two related tasks
carried out by the brain (cf. Dayan and Abbott, \citeyear{DayAbb05}) through two kinds
of models.
The first kind, the encoding model,
predicts brain signals from visual stimuli, while the second kind,
the decoding model recovers visual stimuli from brain signals. Often,
decoding models
are built upon encoding models and hence indirectly validate the
former, but
they are important in their own right. In the September issue of
\emph{Current Biology}, our joint paper with the Gallant Lab, Nishimoto
et al.
(\citeyear{NisVuNas11}) invents a decoding (or movie reconstruction) algorithm
to reconstruct movies from fMRI brain signals. This work has
received intensive and extensive coverage by the media including The Economist's
Oct. 29th 2011 issue (``Reading the Brain: Mind-Goggling'') and
the National Public Radio in their program ``Forum
with Michael Krasny'' on Tue, Sept. 27, 2011 at 9:30 am (``Reconstructing
the Mind's Eye'').
As mentioned earlier, it was selected by the Time Magazine as one of
the best 50 inventions of 2011 and dubbed as ``Mind-reading Computers'' on
the cover page of the Time's invention issue.\vspace*{6pt}

\emph{What is really behind the movie reconstruction algorithm?}\vspace*{6pt}

\emph{Can we learn something from it about how brain works?}\vspace*{6pt}

The movie reconstruction algorithm consists of
statistical encoding and decoding models, both of which
employ regularization. The former are sparse models so they are concise
enough to be viewed and are built
via Lasso${}+{}$CV for each voxel separately. However, as is well-known
Lasso${}+{}$CV results are not stable or reliable enough for scientific
interpretation due to
the $L_1$ regularization and the emphasis of CV on prediction performance.
So Lasso${}+{}$CV is not estimation stable.
The decoding model uses the estimated encoding model for each voxel
and Tiknohov regularization or Ridge in covariance estimation to
pull information across different voxels over V1, V2 and V4 (Nishimoto
et al., \citeyear{NisVuNas11}).
Then an empirical prior for clips of short videos is used
from movie trailers and YouTube to induce posterior weights
on video clips in the empirical prior database.
Tiknohov or Ridge regularization concerns itself with the estimation
of the covariance between voxels that is not of interest
for interpretation. The encoding phase is the focus here from now
on.\looseness=-1

V1 is a primary visual cortex area and the best understood area
in the visual cortex.
Hubel and Wiesel received a Nobel Prize in Physiology or Medicine in 1981
for two major scientific discoveries. One is Hubel and Wiesel (\citeyear{HubWie59})
that uses cat physiology data to show, roughly speaking, that simple V1
neuron cells act like Gabor filters or as angled edge detectors.
Later, using solely image data, Olshausen and Field (\citeyear{OlsFie96}) showed that
image patches can be sparsely represented on Gabor-like basis image patches.
The appearance of Gabor filters in both places is likely not a coincidence,
due to the fact that our brain has evolved to represent the natural world.
These Gabor filters have different locations, frequencies and orientations.
Previous work from the Gallant Lab has built a filter-bank of such Gabor
filters and successfully used them to design encoding models
with single neuron signals in V1 invoked by static natural image stimuli
(Kay et al., \citeyear{KayNasPre08}, Naselaris et al., \citeyear{NasKayGal11}).

In Nishimoto et al. (\citeyear{NisVuNas11}), we use fMRI brain signals observed
over 2700 voxels in different areas of the visual cortex.
fMRI signals are indirect and non-invasive measures of neural activities
in the brain and have good spatial coverage and temporal
resolution in seconds. Each voxel is roughly a cube of 1~mm by 1~mm by
1~mm
and contains hundreds of thousands of neurons. Leveraging the success
of Gabor-filter based models for single neuron brain signals, for a
given image,
a vector of features is extracted by 2-d wavelet filters.
This feature vector has been used to build encoding models for fMRI
brain signals in Kay et al. (\citeyear{KayNasPre08}) and Naselaris et al. (\citeyear{NasKayGal11}).
Invoked by clips of videos/movies, fMRI signals from three subjects are
collected with the same experimental set-up.
To model fMRI signals invoked by movies, a 3-dim motion-energy Gabor
filter bank
has been built in Nishimoto et al. (\citeyear{NisVuNas11}) to extract a feature vector
of dimension of 26K.
Linear models are then built on these features at the observed time
point and
lagged time points.

At present sparse linear regression models are favorites of the Gallant
Lab through Lasso or
$\varepsilon$-L2Boost. These sparse models give similar prediction
performance on validation data as neural nets and kernel machines
on image-fMRI data; they correspond well
to the neuroscience knowledge on V1; and they are easier to interpret than
neural net and kernel machine models that include all features or variables.

For each subject, following a rigorous protocol in the Gallant Lab, the
movie data
(how many frames per second?)
consists of three batches: training, test and validation.
The training data is used to fit a sparse encoding model
via Lasso or e-L2Boost and the test data is used to select
the smoothing parameter by CV. These data are averages of two or three
replicates.
That is, the same movie is played to one subject two or three times and
the resulted fMRI signals are called replicates.
Then a completed encoding determined model is used to predict
the fMRI signals in the validation data (with 10$+$ replicates) and
the prediction performance is measured by the correlation
between the predicted fMRI signals and observed fMRI signals,
for each voxel and for each subject. Good prediction performance
is observed for such encoding models (cf. Figure~\ref{Fig:scatter-hist}).

\section{Stability considerations in the literature}\label{Sec:stability}

Prediction and movie reconstruction are good steps to validate the
encoding model in order to understand the human visual pathway. But the
science lies in finding the features that might drive a voxel, or to
use Tukey's words, finding out ``how things really are.''

It is often the case that the number of data units is easily different
from what is in collected data. There are some hard resource
constraints such as that human subjects can not lie inside an fMRI
machine for too long and it also costs money to use the fMRI machine.
But whether the data collected is for 2 hours as in the data or 1 hours
50 min or 2 hours and 10 min is a judgement call by the experimenter
given the constraints. Consequently, scientific conclusions, or in our
case, candidates for driving features, should be stable relative to
removing a small proportion of data units, which is one form of
reasonable or appropriate data perturbation, or reproducible without a
small proportion of the data units. With a smaller set of data, a more
conservative scientific conclusion is often reached, which is deemed
worthwhile for the sake of more reliable results.

Statistics is not the only field that uses mathematics to describe
phenomena in the natural world. Other such fields include numerical
analysis, dynamical systems and PDE and ODE. Concepts of stability are
central in all of them, implying the importance of stability in
quantitative methods or models when applied to real world problems.

The necessity for a procedure to be robust to data perturbation is a
very natural idea, easily explainable to a child. Data perturbation has
had a long history in statistics, and it has at least three main forms:
jacknife, sub-sampling and bootstrap. Huber (\citeyear{Hub02}) writes in ``John W.
Tukeys Contribution to Robust Statistics:

``[Tukey] preferred to rely on the actual batch of data at hand rather
than on a hypothetical underlying population of which it might be a
sample'' (p. 1643).

All three main forms of data perturbation rely on an ``actual batch of
data'' even though their theoretical analyses do assume hypothetical
underlying populations of which data is a sample. They all have had
long histories.

Jacknife can be traced back at least to Quenouille (\citeyear{Que49,Que56}) where
jacknife was used to estimate the bias of an estimator. Tukey (\citeyear{Tuk58}),
an abstract in the Annals of Mathematical Statistics, has been regarded
as a key development because of his use of jacknife for variance
estimation. Miller (\citeyear{Mil74}) is an excellent early review on Jacknife with
extensions to regression and time series situations. Hinkley (\citeyear{Hin77})
proposes weighted jacknife for unbalanced data for which Wu (\citeyear{Wu86})
provides a theoretical study. K\"{u}nsch (\citeyear{Kun89}) develops Jacknife
further for time series. Sub-sampling on the other hand was started
three years earlier than jacknife by Mahalanobis (\citeyear{Mah46}). Hartigan
(\citeyear{Har69,Har75}) buids a framework for confidence interval estimation based
on subsampling. Carlstein (\citeyear{Car86}) applies subsampling (which he called
subseries) to the time series context. Politis and Romano (\citeyear{PolRom92}) study
subsampling for general weakly dependent processes. Cross-validation
(CV) has a more recent start in Allen (\citeyear{All74}) and Stone (\citeyear{Sto74}). It gives
an estimated prediction error that can be used to select a particular
model in a class of models or along a path of regularized models. It
has been wildly popular for modern data problems, especially for
high-dimensional data and machine learning methods. Hall (\citeyear{Hal83}) and Li
(\citeyear{Li86}) are examples of theoretical analyses of CV. Efron's (\citeyear{Efr79})
bootstrap is widely used and it can be viewed as simplified jacknife or
subsampling. Examples of early theoretical studies of bootstrap are
Bickel and Freedman (\citeyear{BicFre81}) and Beran (\citeyear{Ber84}) for the i.i.d. case, and K\"{u}nsch (\citeyear{Kun89}) for time series. Much more on these three data
perturbation schemes can be found in books, for example, by Efron and
Tibshirani (\citeyear{EfrTib93}), Shao and Tu (\citeyear{ShaTu95}) and Politis, Romano and Wolf (\citeyear{PolRomWol99}).

If we look into the literature of probability theory, the mathematical
foundation of statistics, we see 5
that a perturbation argument is central to limiting law results such as
the Central Limit Theorem (CLT).

The CLT has been the bedrock for classical statistical theory. One
proof of the CLT that is composed of two steps and is well exposited in
Terence Tao's lecture notes available at his website (Tao, \citeyear{Tao12}). Given
a normalized sum of i.i.d. random variables, the first step proves the
universality of a limiting law through a perturbation argument or the
Lindebergs swapping trick. That is, one proves that a perturbation in
the (normalized) sum by a random variable with matching first and
second moments does not change the
(normalized) sum distribution. The second step finds the limit law by
way of solving an ODE.\looseness=-1

Recent generalizations to obtain other universal limiting distributions
can be found in Chatterjee (\citeyear{Cha06}) for Wigner law under non-Gaussian
assumptions and in Suidan (\citeyear{Sui06}) for last passage percolation. It is
not hard to see that the cornerstone of theoretical high-dimensional
statistics, concentration results, also assumes stability-type
conditions. In learning theory, stability is closely related to good
generalization performance (Devroye and Wagner, \citeyear{DevWag79}, Kearns and Ron,
\citeyear{KeaRon99}, Bousquet and Elisseeff, \citeyear{BouEli02}, Kutin and
Niyogi, \citeyear{KutNiy02}, Mukherjee et al., \citeyear{MukNiyPogRif06}, Shalev-Shwartz et al., \citeyear{ShaShaSre10}).

To further our discussion on stability, we would like to explain what
we mean by statistical stability. We
say statistical stability holds if statistical conclusions are robust
or stable to appropriate perturbations to data. That is, statistical
stability is well defined relative to a particular aim and a particular
perturbation to data (or model). For example, aim could be estimation,
prediction or limiting law. It is not difficult to have statisticians
to agree on what are appropriate data perturbations when data units are
i.i.d. or exchangeable in general, in which case subsampling or
bootstrap are appropriate. When data units are dependent,
transformations of the original data are necessary to arrive at
modified data that are close to i.i.d. or exchangeable, such as in
parametric bootstrap in linear models or block-bootstrap in time
series. When subsampling is carried out, the reduced sample size in the
subsample does have an effect on the detectable difference, say between
treatment and control. If the difference size is large, this reduction
on sample size would be negligible. When the difference size is small,
we might not detect the difference with a reduced sample size, leading
to a more conservative scientific conclusion. Because of the utmost
importance of reproducibility for science, I believe that this
conservatism is acceptable and may even be desirable in the current
scientific environment of over-claims.

\section{Estimation stability: Seeking more stable models than Lasso${}+{}$CV}
\label{Sec:ES-CV}

For the fMRI problem, let us recall that
for each voxel, Lasso or e-L2Boost is used to fit the mean function in
the encoding
model with CV to choose the smoothing parameter.
Different model selection criteria have been known to be unstable.
Breiman (\citeyear{Bre96}) compares predictive stability among forward
selection, two versions of garotte and Ridge and their stability
increases in that order. He goes on to propose averaging
unstable estimators over different perturbed data sets in
order to stabilize unstable estimators. Such estimators are
prediction driven, however, and
they are not sparse and thereby not suitable for interpretation.

In place of bootstrap for prediction error estimation
as in Efron (\citeyear{Efr82}), Zhang (\citeyear{Zha93}) uses multi-fold cross-validation
while Shao (\citeyear{Sha96}) uses m out of n bootstrap samples with $m\ll n$.
They then select models with this estimated prediction error, and
provide theoretical results for low dimensional
or fixed p linear models. Heuristically,
the m out of n bootstrap in Shao (\citeyear{Sha96}) is needed because the model selection
procedure is a discrete (or set) valued estimator for the true
model predictor set and hence non-smooth (cf.~Bickel, G\"otze, and van
Zwet, \citeyear{Bickeletal97}).

The Lasso (Tibshirani, \citeyear{Tib96}) is a modern model selection method for
linear regression and very popular in high-dimensions:
\[
\hat{\beta} (\lambda) = \mathop{\arg} _{\beta\in R^p} \bigl\{\Vert Y-X \beta\Vert _2^2
+\lambda\Vert \beta\Vert _1 \bigr\},
\]
where $Y \in R^n$ is the response vector and $X \in R^{n\times p}$
is the design matrix. That is, there are $n$ data units and
$p$ predictors. For each $\lambda$, there is a unique $L_1$ norm
for its solution that we can use to index the solution as
$\hat{\beta} (\tau)$ where
\[
\tau= \tau(\lambda) = \bigl\Vert \hat{\beta} (\lambda)\bigr\Vert _1.
\]

Cross-validation (CV) is used most of the time to select $\lambda$ or
$\tau$, but
Lasso${}+{}$CV is unstable relative to bootstrap
or subsampling perturbations when predictors are correlated (cf.
Meinshausen and B\"{u}hlmann, \citeyear{MeiBul10}, Bach, \citeyear{Bac08}).


Using bootstrap in a different manner than Shao (\citeyear{Sha96}),
Bach (\citeyear{Bac08}) proposes BoLasso to improve
Lasso's model selection consistency property by taking
the smallest intersecting model of selected models over
different bootstrap samples. For particular smoothing parameter sequences,
the BoLasso selector is shown by Bach (\citeyear{Bac08}) to be model selection consistent
for the low dimensional case without the irrepresentable
condition needed for Lasso (cf. Meinshausen and B\"{u}hlmann, \citeyear{MeiBul06},
Zhao and Yu, \citeyear{ZhaYu06}; Zou, \citeyear{Zou06}; Wainwright, \citeyear{Wai09}).
Meinshausen and B\"{u}hlmann (\citeyear{MeiBul10}) also weaken the irrepresentable condition
for model selection consistency of a stability selection criterion built
on top of Lasso. They bring perturbations to a Lasso path
through a random scalar vector in the Lasso $L_1$ penalty, resulting
in many random Lasso paths. A threshold parameter is needed
to distinguish important features based on these random paths. They do
not consider the problem of selecting one smoothing parameter value for Lasso
as in Lim and Yu (\citeyear{LimYu13}).

We would like to seek a specific model along the Lasso path to
interpret and hence
selects a specific $\lambda$ or $\tau$. It is well known that CV does
not provide
a good interpretable model because Lasso${}+{}$CV is unstable. Lim and Yu
(\citeyear{LimYu13}) propose a stability-based
criterion that is termed Estimation Stability (ES).
They use the cross-validation data perturbation scheme. That is,
$n$ data units are randomly partitioned into
$V$ blocks of pseudo data sets of size $(n-d)$ or subsamples
where $d = \lfloor n/V \rfloor.$\footnote{$\lfloor x \rfloor$ is the
floor function or
the largest integer that is smaller than or equal to $x$.}

Given a smoothing parameter $\lambda$, a Lasso estimate $\hat{\beta}_v
(\lambda)$
is obtained for the $v$th block $v=1,\ldots,V$.
Since the $L_1$ norm is a meaningful quantity to line up the $V$
different estimates, Lim and Yu (\citeyear{LimYu13})\footnote{It is also fine to use
$\lambda$
to line up the different solutions, but not a good idea to use the
ratio of $\lambda$
and its maximum value for each pseudo data set.} use it, denoted as
$\tau$ below, to line
up these estimates
to form an estimate $\hat{m} (\tau)$
for the mean regression function and an approximate delete-$d$
jacknife estimator for the variance of $\hat{m} (\tau)$:
\begin{eqnarray*}
\hat{m} (\tau)&=& \frac{1}{V} \sum_v X \hat{
\beta}_v (\tau),
\\
\hat{T}(\tau) &=&\frac{n-d}{d} \frac{1}{V} \sum
_v \bigl(\bigl\Vert  X \hat{\beta}_v (\tau) - \hat{m} (
\tau)\bigr\Vert ^2\bigr).
\end{eqnarray*}
The last expression is only an approximate delete-$d$ jacknife variance
estimator unless $V={n \choose n-d}$ when all the subsamples of size $n-d$
are used.
Define the (estimation) statistical stability measure as
\[
ES(\tau) = \frac{{1}/{V} \sum_v \Vert  X \hat{\beta}_v (\tau) -
\hat{m} (\tau)\Vert ^2 }{ \hat{m}^2 (\tau)} = \frac{d}{n-d}\frac{\hat{T}(\tau)}{\hat{m}^2(\tau)} =
\frac{d}{n-d} \frac{1}{Z^2 (\tau)} ,
\]
where $Z (\tau)=\hat{m}(\tau)/\sqrt{\hat{T}(\tau)}$.

For nonlinear regression functions, ES can still be applied
if we take an average of the estimated regression functions.
Note that ES aims at estimation stability, while CV aims at prediction
stability.
In fact, ES is the reciprocal of a test statistic for testing
\[
H_0\dvtx X \beta= 0.
\]
Since $Z (\tau)=\hat{m}(\tau)/\sqrt{\hat{T}(\tau)}$ is a test statistic
for $H_0$,
$Z^2 (\tau)$ is also a test statistic.
$ES(\tau)$ is a scaled version of the reciprocal $1/Z^2 (\tau)$.

To combat the high noise situation where ES would not
have a well-defined minimum, Lim and Yu (\citeyear{LimYu13}) combine ES with CV to
propose the
\textit{ES-CV selection criterion} for smoothing parameter $\tau$:\vspace*{6pt}

\emph{Choose the largest $\tau$ that minimizes ES $(\tau)$ and is
smaller or equal
to the CV selection.}\vspace*{6pt}

\begin{figure}

\includegraphics{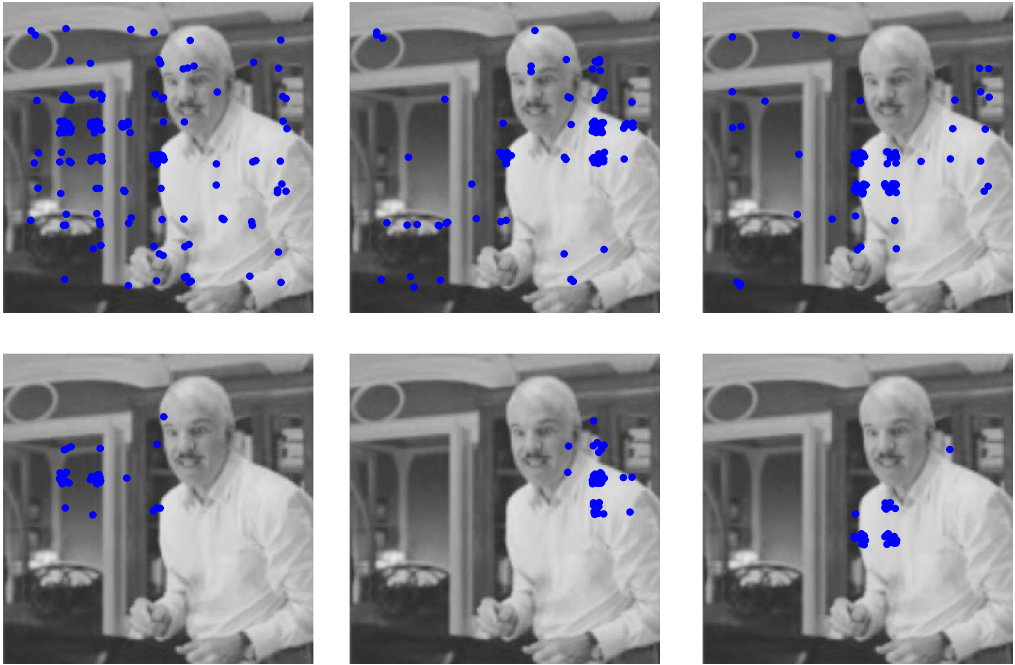}

\caption{For three voxels (one particular subject), we display the
(jittered)
locations that index the Gabor features selected by
CV-Lasso (top row) and ESCV-Lasso (bottom row).}
\label{Fig:feature-locations}
\end{figure}

ES-CV is applicable to smoothing parameter selection in Lasso,
and other regularization methods such as
Tikhonov or Ridge regularization (see, for example, Tikhonov, \citeyear{Tik43},
Markovich, \citeyear{Mar07}, Hoerl, \citeyear{Hoe62}, Hoerl and Kennard, \citeyear{HoeKen70}). ES-CV is well suited
for parallel computation as CV and incurs only a negligible computation
overhead because $\hat{m} (\tau)$ are already computed for CV.
Moreover, simulation studies in Lim and Yu (\citeyear{LimYu13}) indicate that,
when compared with Lasso${}+{}$CV,
ES-CV applied to Lasso gains dramatically in terms of false discovery rate
while it loses only somewhat in terms of true discovery rate.

The features or predictors in the movie-fMRI problem
are 3-d Gabor wavelet filters, and each of them is characterized
by a (discretized) spatial location on the image, a (discretized) frequency
of the filter, a (discretized)
orientation of the filter, and 4 (discrete) time-lags on the corresponding
image that the 2-d filter is acting on.
For the results comparing CV and ES-CV in Figure~\ref{Fig:feature-locations},
we have a sample size $n=7\mbox{,}200$ and use a reduced set of $p = 8\mbox{,}556$ features
or predictors,
corresponding to a coarser set of filter frequencies than what is
used in Nishimoto et al. (\citeyear{NisVuNas11}) with $p=26\mbox{,}220$ predictors.

\begin{figure}

\includegraphics{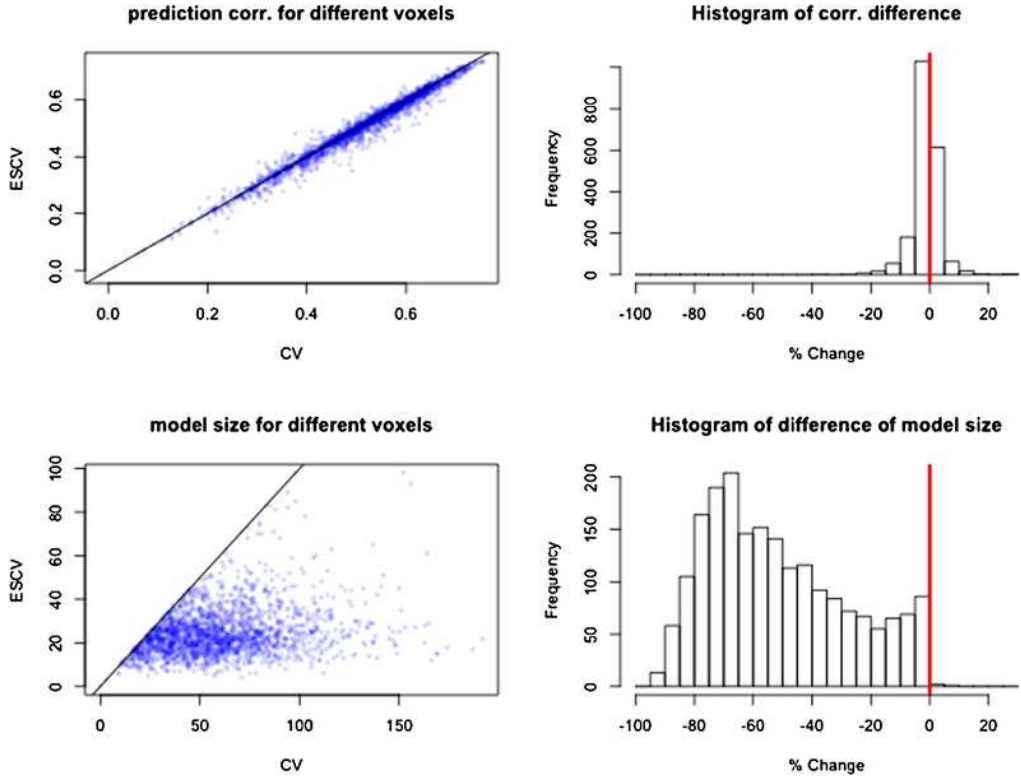}

\caption{Comparisons of ESCV(Lasso) and CV(Lasso) in terms of model
size
and prediction correlation. The scatter plots on the left compare ESCV
and CV
while the histograms on the right display the differences of model size
and prediction correlation.}\vspace*{-2pt}
\label{Fig:scatter-hist}
\end{figure}

We apply both CV and ES-CV to select the smoothing parameters
in Lasso (or e-L2Boost).
For three voxels (and a particular subject),
for the simplicity of display,
we show the locations of the selected
features (regardless of their frequencies, orientations and
time-lags) in Figure~\ref{Fig:feature-locations}. For these three voxels,
ES-CV maintains
almost the same prediction correlation performances as CV (0.70 vs.
0.72) while
ES-CV selects many fewer and more concentrated locations than CV.
Figure~\ref{Fig:scatter-hist} shows the comparison results
across 2088 voxels in the visual cortex that are selected for their high
SNRs. It is composed of four sub-plots. The upper two plots compare prediction
correlation performance of the models built via Lasso with CV and ES-CV
on validation data.
For each model fitted on training data and each voxel, predicted responses
over the validation data are calculated. Its correlation with the observed
response vector is the ``prediction correlation'' displayed in Figure \ref
{Fig:scatter-hist}. The lower two plots
compare the sparsity properties of the models or model size.
Because of the definition of ES-CV, it is expected that the ES-CV model
are always smaller than or equal to the CV model. The sparsity advantage
of ES-CV is apparent with a huge overall reduction of 60\%
on the number of selected features and a minimum loss of overall
prediction accuracy
by only 1.3\%. The average size of the ES-CV models is 24.3 predictors,
while that for the CV models is 58.8 predictors; the average prediction
correlation performance of the ES-CV models is 0.499, while that for
the CV models is 0.506.


%
%

\section{Sample variability meets robust statistics in high-dimensions}
\label{Sec:robust}

Robust statistics also deals with stability, relative to model perturbation.
In the preface of his book ``Robust Statistics,'' Huber (\citeyear{Hub81}) states:

``Primarily, we are concerned with \emph{distributional robustness}:
the shape
of the true underlying distribution deviates slightly from the assumed model.''

Hampel, Rousseeuw, Ronchetti and Stahel (\citeyear{HamRonRou86}) write:

``Overall, and in analogy with, for example, the stability aspects of
differential equations or of numerical computations, robustness
theories can be viewed as stability theories of statistical inference''
(p. 8).

Tukey (\citeyear{Tuk58}) has generally been regarded as the first paper
on robust statistics. Fundamental contributions were made by
Huber (\citeyear{Hub64}) on M-estimation of location parameters, Hampel (\citeyear{Ham68,Ham71,Ham74}) on
``break-down'' point and influence curve. Further important contributions
can be found, for example, in Andrews et al. (\citeyear{Andrewsetal72}) and Bickel (\citeyear{Bic75})
on one-step Huber estimator,
and in Portnoy (\citeyear{Por77}) for M-estimation in the dependent
case.\looseness=-1

For most statisticians, robust statistics in linear regression is
associated with studying estimation problems when the errors have
heavier tail distributions than
the Gaussian distribution. In the fMRI problem, we fit mean functions
with an L2 loss.
What if the ``errors'' have heavier tails than
Gaussian tails? For the L1 loss is commonly used in robust statistics
to deal with heavier tail errors in regression, we may wonder
whether the L1 loss would add more stability to the fMRI problem.
In fact, for high-dimensional data such as in our fMRI problem,
removing some data units could severely change the outcomes of our model
because of feature dependence.
This phenomenon is also seen in simulated data from linear models with
Gaussian errors in high-dimensions.

\emph{How does sample to sample variability interact with heavy tail
errors in high-dimensions?}

In our recent work El Karoui et al. (\citeyear{ElKarouietal12}), we
seek insights into this question through analytical work.
We are able to see interactions between sample variability and
double-exponential tail errors
in a high-dimensional linear regression model. That is,
let us assume the following linear regression model
\[
Y_{n\times1} = X_{n \times p} \beta_{p \times1} + \varepsilon_{n \times1},
\]
where
\[
X_i \sim N(0, \Sigma_p), \mbox{i.i.d.} ,
\varepsilon_i \mbox{i.i.d.}, E\varepsilon_i =0, E
\varepsilon_i^2 = \sigma^2 < \infty.
\]

An M-estimator with respect to loss function $\rho$ is given as
\[
\hat{\beta} = \mathop{\operatorname{argmin}}_{\beta\in R^p} \sum_i \rho
\bigl(Y_i- X_i^\prime\beta\bigr).
\]

We consider the random-matrix high-dimensional regime:
\[
p/n \rightarrow\kappa\in(0,1).
\]

Due to rotation invariance, WLOG, we can assume $\Sigma_p = I_p$ and
$\beta=0$. We cite
below a result from El Karoui et al. (\citeyear{ElKarouietal12}) for the important
special case of $\Sigma_p = I_p$:

\begin{result}[(El Karoui et al., \citeyear{ElKarouietal12})]\label{res1}
Under the aforementioned assumptions, let $r_\rho(p, n) = \Vert \hat{\beta
}\Vert $, then $\hat{\beta}$ is distributed as
\[
r_\rho(p, n) U,
\]
where $U \sim \operatorname{uniform} (S^{p-1}) (1)$, and
\[
r_\rho(p, n) \rightarrow r_\rho(\kappa) ,
\]
as $n, p \rightarrow\infty$ and $p/n
\rightarrow\kappa\in(0,1)$.\vadjust{\goodbreak}
\end{result}

Denote
\[
\hat{z}_{\varepsilon} := \varepsilon+ r_\rho(\kappa) Z,
\]
where $Z \sim N(0,1)$ and independent of $\varepsilon$, and let
\[
\operatorname{prox}_c (\rho) (x) = \mathop{\operatorname{argmin}}_{ y \in R} \biggl[ \rho(y) +
\frac{(x-y)^2}{2c} \biggr] .
\]

Then $r_\rho(\kappa)$ satisfies the following system of equations
together with some nonnegative~$c$:
\begin{eqnarray*}
E \bigl\{ \bigl[ \operatorname{prox}_c (\rho)\bigr]^\prime\bigr\}& =& 1 -
\kappa,
\\
E \bigl\{ \bigl[ \hat{z}_\varepsilon- \operatorname{prox}_c (
\hat{z}_\varepsilon) \bigr]^2\bigr\} &=& \kappa r_\rho^2
(\kappa) .
\end{eqnarray*}

In our limiting result, the norm of an M-estimator stabilizes.
It is most interesting to mention that in the proof
a ``leave-one-out'' trick is
used both row-wise and column-wise such that
one by one rows are deleted and similarly columns are deleted.
The estimators with deletions are then compared to the estimator
with no deletion.
This is in effect a perturbation argument and reminiscent of the
``swapping trick''
for proving the CLT as discussed before. Our
analytical derivations involve prox functions, which
are reminiscent of the second step in proving normality in the CLT.
This is because a prox function is a form of derivative,
and not dissimilar to the derivative appearing in the ODE derivation
of the analytical form of the limiting distribution (e.g normal
distribution) in the CLT.

In the case of i.i.d. double-exponential errors, El Karoui et al. (\citeyear{ElKarouietal12})
numerically solve the two equations in Result~\ref{res1} to show that
when $\kappa\in(0.3,1)$,
$L_2$ loss fitting (OLS) is better than $L_1$ loss fitting (LAD) in
terms of MSE or variance. They also show that the numerical results
match very well
with simulation or Monte Carlo results. At a high level,
we may view that
$\hat{z}_\varepsilon$ holds the key
to this interesting phenomenon.
Being a weighted convolution of $Z$ and $\varepsilon$, it
embeds the interaction
between sample variability (expressed in $Z$) and error variability
(expressed in $\varepsilon$) and this interaction is captured
in the optimal loss function (cf. El Karoui et al., \citeyear{ElKarouietal12}).
In other words,
$\hat{z}_\varepsilon$ acts more like double exponential when the influence
of standard normal $Z$ in $\hat{z}_\varepsilon$ is not dominant (or when
$\kappa< 0.3$ or so as we discover
when we solve the equations) and in this case, the optimal loss
function is closer to LAD loss.
In cases when $\kappa> 0.3$,
it acts more like Gaussian noise, leading
to the better performance of OLS (because the optimal loss is closer
to LS).


Moreover, for double exponential errors,
the M-estimator LAD is an MLE and we are in a high-dimensional situation. It is well-known that MLE does not work in
high-dimensions.
Remedies have been found through penalized MLE where a bias is introduced
to reduce variance and consequently reduce the MSE.
In contrast, when $\kappa\in(0.3, 1)$, the better estimator OLS is
also unbiased, but has a smaller variance nevertheless. The variance reduction
is achieved through a better loss function LS than the LAD and
because of a concentration of quadratic forms
of the design matrix. This concentration does not hold
for fixed orthogonal designs, however.
A~follow-up work (Bean et al., \citeyear{Beanetal12}) addresses the question of
obtaining the optimal loss function. It is current research regarding
the performance of estimators from penalized OLS and penalized LAD
when the error distribution is double-exponential.
Preliminary results indicate that similar phenomena occur
in non-sparse cases.

Furthermore, simulations with design matrix from an fMRI experiment and
double-exponential
error show the same phenomenon, that is, when $\kappa=p/n > 0.3$ or so,
OLS is better than LAD. This provides some insurance for using L2 loss
function in the fMRI project.
It is worth noting that El Karoui et al. (\citeyear{ElKarouietal12}) contains results for
more general settings.
\section{Conclusions}
\label{Sec:conclusions}

In this paper, we cover three problems facing statisticians at the 21st
century: figuring out how vision works with fMRI data,
developing a smoothing parameter selection method for Lasso, and
connecting perturbation in the case of high-dimensional data with classical
robust statistics through analytical work.
These three problems are tied together by stability.
Stability is well defined if we describe the data perturbation
scheme for which stability is desirable, and such schemes include bootstrap,
subsampling, and cross-validation. Moreover, we briefly review
results in the probability literature to explain that stability is driving
limiting results such as the Central Limit Theorem, which is
a foundation for classical asymptotic statistics.

Using these three problems as backdrop, we make four points.
Firstly, statistical stability considerations can effectively aid
the pursuit for interpretable and reliable scientific models, especially
in high-dimensions. Stability in a broad sense includes replication,
repeatability, and
different data perturbation schemes.
Secondly, stability is a general principle on which to build
statistical methods for different purposes.
Thirdly, the meaning of stability needs articulation in high-dimensions
because it could be brought about by sample variability and/or heavy tails
in the errors of a linear regression model.
Last but not least, emphasis should be placed on the stability aspects
of statistical inference and conclusions, in the referee process of
scientific and applied statistics papers and in our current statistics
curriculum.

Statistical stability in the age of massive data is an important area for
research and action because high-dimensions provide ample opportunities
for instability
to reveal itself to challenge reproducibility of scientific findings.

As we began this article with words of Tukey, it seems fitting to end
also with his words:

\begin{quote}
``What of the future? The future of data analysis can involve great
progress, the
overcoming of real difficulties, and the provision of a great service
to all fields of science and technology. Will it? That remains to us,
to our willingness to take up the
rocky road of real problems in preferences to the smooth road of unreal
assumptions, arbitrary
criteria, and abstract results without real attachments. Who is for the
challenge?'' -- Tukey (p. 64, \citeyear{Tuk62}).
\end{quote}


\section*{Acknowledgements}
This paper is based on the 2012 Tukey Lecture of the Bernoulli Society
delivered by the author at the 8th
World Congress of Probability and Statistics in Istanbul on July 9, 2012.
For their scientific influence and friendship,
the author is indebted to her teachers/mentors/colleagues,
the late Professor Lucien Le Cam, Professor Terry Speed, the late Professor
Leo Breiman, Professor Peter Bickel, and Professor Peter B\"uhlmann.
This paper is invited for the special issue of \emph{Bernoulli}
commemorating the 300th
anniversary of the publication of Jakob Bernoulli's Ars Conjectandi in 1712.

The author would like to thank Yuval Benjamini for his help on
generating the results in the figures. She would also like to thank two
referees for their detailed and insightful comments, and Yoav Benjamini
and Victoria Stodden for helpful discussions. Partial supports are
gratefully acknowledged by NSF Grants SES-0835531 (CDI) and
DMS-11-07000, ARO Grant W911NF-11-1-0114, and the NSF Science and
Technology Center on Science of Information through Grant CCF-0939370.


%



\printhistory

\end{document}